\documentclass[12pt,english]{article}
\usepackage[T1]{fontenc}
\usepackage{a4wide}
\usepackage{graphicx}
\usepackage{setspace}
\usepackage{amssymb}

\makeatletter

\providecommand{\LyX}{L\kern-.1667em\lower.25em\hbox{Y}\kern-.125emX\@}
 \usepackage{verbatim}

\usepackage{amsthm}

\newcommand{\ds}{\displaystyle}

\newtheorem{theorem}{Theorem}

\newtheorem{lemma}[theorem]{Lemma}
\newtheorem{cor}[theorem]{Corollary}
\newtheorem{rem}[theorem]{Remark}
\newtheorem{conj}[theorem]{Conjecture}

\usepackage{babel}
\makeatother
\begin{document}

\newcommand{\bZ}{\mathbf{Z}}
 
\newcommand{\Z}{\bZ }
 
\newcommand{\cM}{\mathcal{M}}
 
\newcommand{\sP}{\mathsf{P}}
 
\newcommand{\Xd}[1]{X(#1)}
 
\newcommand{\Yd}[1]{Y(#1)}
 
\newcommand{\sE}{\mathsf{E}\, }
 
\newcommand{\sgeq}{\geq _{{st}}}
 
\newcommand{\sleq}{\leq _{\mathrm{s}t}}
 
\newcommand{\R}{\mathbf{R}}
 
\newcommand{\be}{\mathbf{e}}
 
\newcommand{\vs}{v^{*}}
 
\newcommand{\cD}{\mathcal{D}}
 
\newcommand{\Xc}{Y}
 
\newcommand{\nconst}{\relax }
 
\newcommand{\flam}{\varphi }
 
\newcommand{\skp}[2]{\left\langle #1,#2\right\rangle }

\newcommand{\ga}[1]{\gamma _{\ref {#1}}}

\newcommand{\CC}[1]{C_{\ref {#1}}}

\newcommand{\IF}[1]{I_{\left\{ \Xc (n+1)\in F_{#1}\right\} }}

\newcommand{\IR}[1]{I_{\left\{ \Xc (n+1)\in R_{#1}\right\} }}

\title{On the cascade rollback synchronization}

\author{Anatoli Manita\thanks{The work of this author is supported by Russian Foundation of Basic Research (grant 02-01-00945). Postal address: Faculty of Mathematics and Mechanics, Moscow State University, 119992, Moscow, Russia. E-mail:~manita@mech.math.msu.su}\and
Fran\c{c}ois  Simonot \thanks{IECN, Universit\'e Henri Poincar\'e Nancy~I, Esstin,  2, Rue  J.~Lamour, 54500 Vandoeuvre, France. E-mail: francois.simonot@esstin.uhp-nancy.fr}
   }

\maketitle
\begin{abstract}
We consider a cascade model of $N$ different processors performing
a distributed parallel simulation. The main goal of the study is to
show that the long-time dynamics of the system has a cluster behavior.
To attack this problem we combine two methods: stochastic comparison
and Foster--Lyapunov functions.
\end{abstract}

\section{Introduction}

The present paper contains a probabilistic analysis of some mathematical
model of asynchronous algorithm for parallel simulation. For the detailed
discussion of synchronization issues in parallel and distributed algorithms
we refer to~\cite{JeWI,BeTsi}. Here we give only a brief description
of the problem. In large-scale parallel computation it is necessary
to coordinate the activities of different processor which are working
together on some common task. Usually such coordination is implemented
by using a so-called message-passing system. This means that a processor
shares data with other processors by sending timestamped messages.
Between sending or receiving the messages the processors work independently.
It can be happened that till the moment of receiving of a message
some processor can proceed farther in performing its program than
the value of timestamp indicated in this newly received message; in
this case the state of the processor should be \emph{rolled back}
to the indicated value. It is clear that due to possible rollbacks
the mean speed of a given processor in the computing network will
be lower than its proper speed. One of the most important performance
characteristics of the system is the progress of the computing network
on large time intervals.

Probabilistic models for such system are studied already for twenty
years. From the probabilistic point of view these models consists
of many relatively independent components which synchronize from time
to time their states according to some special algorithm. The detailed
review of all existing publications is out of range of this paper.
We would like to mention only that the bibliography on this subject
consists mostly of two group of papers. The first group of publication
\cite{MitMit,KuSh,GrShSt,GuKuSh,ShKuRe,GuKu}  are devoted to the
case of two processors. The paper~\cite{MitMit} is of special interest
since it contains an exhaustive basic analysis of the two-dimensional
model and had a big influence on further research. The case of many
processors was studied in \cite{MadWalMes,GuAkFu,AkChFuSer,PGM,GrShSt,T1,T3}.
An important difference of such models from two-dimensional case is
that in a realistic model with more than two processors one message
can provoke a multiple rollback of a chain of processors. Since the
multi-dimensional model is much more complicated for a rigorous study,
in the above papers the authors deal with the set of identical processors
and their mathematical results are contained in preparatory sections
before large numerical simulations. 

It should be noted also that probabilistic models with synchronization
mechanism are interesting also for modelling database systems~(see
for example, \cite{JeWI}). Moreover, now synchronization-like interactions
are considered as well in the framework of interaction particle systems~\cite{MM-RR,MShCh,MM-TVP}. 

The model considered in the present paper is of special interest for
the following reasons. We deals with a nonhomogeneus model consisting
of several different processors. We consider case of message-passing
topology other from the topology of complete graph which was considered
in all previous papers. Our main interest is the cascade model which
pressupose a subordination between processors. We put forward a conjecture
on the cluster behavior of the system: processors can be divided into
separated groups which are asymptotically independent and have their
own proper performance characteristics. Our main goal is to justify
this conjecture. One should point out that in the case of complete
graph topology the cluster decomposition into groups is degenerated
and, thus, not interesting. 

We describe our model in terms of multi-dimensional continuous time
Markov process. To get asymptotical performance characteristics of
the model we combine two probabilistic methods (stochastic comparison
and Foster--Lyapunov functions). 

The paper is organized as follows. In Section~\ref{sec:Description-of-continuous}
we introduce a general continuous time Markov model and define a cascade
model as a special subclass of the general model. In~Section~\ref{sec:Definition-of-discrete}
we pass to the embedded Markov chain. Main problem is now to study
a long-time behavior of Markov chain with highly nonhomogeneous transition
probabilities. To do this we consider relative coordinates and find
groups of processors whose evolution is ergodic (convergences to a
steady state) in these relative coordinates. To our opinion the method
of Foster-Lyapunov functions seems to be the only one to prove the
stability in relative coordinates for the Markov chain under consideration.
First of all in Section~\ref{sec:Case-N-2} we start from the case
of two processors ($N=2$) and the analysis here is rather simple
and similar to~\cite{MitMit}. In the study of the three-dimensional
case (Section~\ref{sec:N-3-Lyap}) the main point is the proof of
ergodicity. We propose an explicit construction of some nonsmooth
Foster-Lyapunov function. Our construction is rather nontrivial as
it can be seen by comparing with already existing explicit examples
of Lyapunov functions (see~\cite{FMM}). All this analysis bring
us to some conclusions presented in Section~\ref{sec:Conclusions}.
This section contains decomposition into groups (clusters) in the
case of cascade model with any number of processors~$N$ and our
main Conjecture~\ref{con-N-main}. We show that the proof of this
conjecture could be related with progress in explicit construction
of multi-dimensional Foster-Lyapunov functions. Analysis of random
walks in~$\Z _{+}^{n}$ (which was done in~\cite{FMM}) shows that,
in general, this technical problem may be very difficult. In the next
papers we hope to overcome these difficulties by using specific features
of our concrete Markov processes.

\paragraph*{Acknowledgements.  }

The first author is very grateful to the team TRIO (INRIA--Lorraine)
and to l'Ecole des Mines de Nancy for their hospitality during his
stay at Nancy in summer 2004 when the main results of this paper were
obtained.

\section{Description of continuous time model}

\label{sec:Description-of-continuous}

\subsection{General model}

We present here some mathematical model for parallel computations.
There are $N$ computing units (processors) working together on some
common task. The state of a processor $k$ is described by an integer
variable $x_{k}\in \bZ $ which is called a local (or inner) time
of the processor~$k$ and has a meaning of amount of job done by
the processor $k$ up to the given time moment. 

Assume that the state $\left(x_{1},x_{2},\ldots ,x_{N}\right)$ of
the system evolves in continuous time~$t\in \R _{+}$. Any change
of a state is possible only at some special random time instants.
Namely, with any processor~$k$ we associate a Poissonian flow $\Pi ^{k}=\left\{ 0=\sigma _{0}^{k}<\sigma _{1}^{k}<\cdots <\sigma _{n}^{k}<\cdots \right\} $
with intensity~$\lambda _{k}$ and with a pair $(k,l)$ of processors
we associate a Poissonian flow $\Pi ^{kl}=\left\{ 0=\sigma _{0}^{kl}<\sigma _{1}^{kl}<\cdots <\sigma _{n}^{kl}<\cdots \right\} $
with intensity~$\beta _{kl}$. This means, for example, that $\left\{ \sigma _{n}^{k}-\sigma _{n-1}^{k}\right\} _{n=1}^{\infty }$
is a sequence of independent exponentially distributed random variables
with mean~$\lambda _{k}^{-1}$: $\forall n=1,2,\ldots \quad \sP \left\{ \sigma _{n}^{k}-\sigma _{n-1}^{k}>s\right\} =\exp \left(-\lambda _{k}s\right)$,
and similarly for the flows $\Pi ^{kl}$. We also assume that all
these flows $\Pi ^{k}$ and $\Pi ^{kl}$ are mutually independent.

Let us now define a stochastic process $\left(X(t)=\left(x_{1}(t),\ldots ,x_{N}(t)\right),t\in \R _{+}\right)$
on the state space $\Z ^{N}$ according to the following rules.

1) At time instants $\sigma _{n}^{k}$ the processor $k$ increases
its local time~$x_{k}$ by $1$:~ $x_{k}(\sigma _{n}^{k}+0)=x_{k}(\sigma _{n}^{k})+1$.

2) There is an exchange of information between different processors.
At time instant~$\sigma _{i}^{kl}$ the processor~$k$ sends a message
$m_{kl}^{(x_{k})}$ to the processor~$l$. We assume that messages
reach their destination immediately. A message $m_{kl}^{(x_{k})}$
coming to node $l$ from node~$k$ contains an information about
local time $x_{k}(\sigma _{i}^{kl})=x_{k}$ of the sender~$k$. \emph{If}
at the time instant~$\sigma _{i}^{kl}$ (when the message $m_{kl}^{(x_{k})}$
arrives to the node~$l$) we have $x_{l}(\sigma _{i}^{kl})>x_{k}(\sigma _{i}^{kl})$
\emph{then} the local time~$x_{l}$ rolls back to the value $x_{k}$:
$x_{l}(\sigma _{i}^{kl}+0)=x_{k}(\sigma _{i}^{kl})$. Moreover, if
the processor~$l$ rolls back, then all messages sent by the processor
$l$ during the time interval~$\mathcal{I}=(\theta _{l}(x_{k},\sigma _{i}^{kl}),\sigma _{i}^{kl}),$
where $\theta _{l}(x,u):=\max \left\{ s\leq u:\, x_{l}(s)=x,\, x_{l}(s+0)=x+1\right\} ,\sigma _{i}^{kl}),$
should be eliminated. This may generate a cascading rollback of local
times for some subset of processors. For example, assume that there
is a processor $q$ which received a message $m_{lq}^{(x'_{l})}$
at some time instant $s'\in \mathcal{I}$ and $x_{q}(\sigma _{i}^{kl})>x_{l}(s')=x'_{l}$.
Then the local clock~of~$q$ should be rolled back to the value
$x_{l}(s')$: $x_{q}(\sigma _{i}^{kl}+0)=x_{l}(s')$ and, moreover,
all messages sent by $q$ during the interval~$\mathcal{I}=(\theta _{q}(x_{l}(s'),\sigma _{i}^{kl}),\sigma _{i}^{kl})$
should be deleted, and so on. Hence, at time instant $\sigma _{i}^{kl}$
a message from~$k$ to $l$ can provoke a multiple rollback of processor~$l,q,\ldots $
in the system.

\subsection{Cascade model}

\label{sub:Cascade-model}

From now we shall consider the following special subclass of the above
general model.

A chain of processors $1,2,\ldots ,N$ is called a \emph{cascade}
if any processor $j$ can send a message only to its right neighbour
$j+1$. Hence, the processor~$N$ does not send any message and the
processor~$1$ does not receive any message. In other words, $\beta _{ij}\not =0\, \Leftrightarrow \, (j=i+1)$.
A message sent from $j$ to $j+1$ can provoke a cascading roll-back
of processors $j+2,\ldots $~. Recall that all above time intervals
are exponentially distributed and assumed to be independent. Obviously,
the stochastic process $X_{c}^{(N)}(t)=\left(\, x_{1}(t),\ldots ,x_{N}(t)\, \right)$
is Markovian. A very important property is that any {}``truncated''
marginal process  $X_{c}^{(N_{1})}(t)=\left(\, x_{1}(t),\ldots ,x_{N_{1}}(t)\, \right)$,
$N_{1}\leq N$, is also Markovian. 

Assume that for any $j$ the following limit \begin{equation}
\vs _{j}=\lim _{t\rightarrow +\infty }\frac{x_{j}(t)}{t}\qquad (\textrm{in probability})\label{eq:v*j-lim}\end{equation}
 exists. Then the numbers $\vs _{j}$, $j=1,\ldots ,N$, characterize
\emph{performance} of the model. The main goal of the present paper
is to prove the existence of these limits and to calculate them.

Note that if we uniformly transform the absolute time scale $t=cs$,
where $c>0$ is a constant and $s$ is a new absolute time scale,
the performance characteristics~(\ref{eq:v*j-lim}) will not change.

\section{Definition of the discrete time cascade model}

\label{sec:Definition-of-discrete}

Consider a sequence \[
0=\tau _{0}<\tau _{1}<\tau _{2}<\cdots <\cdots \]
of time moments when changes of local time at nodes may happen (we
mean local time updates and moments of sending of messages). It is
clear that $\left\{ \tau _{r+1}-\tau _{r}\right\} _{r=0}^{\infty }$
is a sequence of independent identically distributed r.v. having exponential
distribution with parameter\[
Z=\sum _{i=1}^{N}\lambda _{i}+\sum _{i=1}^{N-1}\beta _{i,i+1}.\]
Observing the continuous time Markov process $(x_{1}(t),\ldots ,x_{N}(t))$
at epochs $\tau _{n}$ we get the so-called embedded discrete time
Markov chain $\left\{ \Xd{n},n=0,1,\ldots \right\} $ with state space
$\bZ _{+}^{N}$. In the sequel we will be interested in the long-time
behaviour of the chain $\left\{ \Xd{n},n=0,1,\ldots \right\} $.

\paragraph*{Transition probabilities.}

In the MC $\left\{ \Xd{n},n=0,1,\ldots \right\} $ there are transitions
produced by the free dynamics and transitions generated by rollbacks.
By the free dynamics we mean updating of local times\[
P\left\{ \Xd{n+1}=x+\be _{j}\, |\, \Xd{n}=x\right\} =\lambda _{j}Z^{-1},\quad j=1,\ldots ,N\, ,\]
where $\be _{j}=(0,\ldots ,0,\begin{array}[t]{c}
 1\\
 ^{j}\end{array},0,\ldots ,0)$. It is easy to see that if a state $x=(x_{1},\ldots ,x_{N})$ is
such that for some $j$ ~ $x_{j}<x_{j+1}$ then a message sent from
$j$ to $j+1$ produces a transition of the following form \begin{equation}
(x_{1},\ldots ,x_{j},x_{j+1},\ldots ,x_{l},x_{l+1},\ldots ,x_{N})\rightarrow (x_{1},\ldots ,x_{j},w_{j+1},\ldots ,w_{l},x_{l+1},\ldots ,x_{N})\label{eq:x-tran}\end{equation}
with probability\begin{equation}
Z^{-1}\beta _{j,j+1}\, \prod _{q=j+1}^{l-1}p(w_{q},x_{q};w_{q+1},x_{q+1})\, \times \, \left(1-b_{l}\right)^{\min \left(x_{l},x_{l+1}-1\right)-w_{l}+1},\label{eq:pr-x-tran}\end{equation}
where

\begin{itemize}
\item sequence $\left(w_{j+1},\ldots ,w_{l}\right)$ is admissible in the
following sense:\[
j<l\leq N,\quad \quad w_{j+1}=x_{j}\qquad w_{q}\leq w_{q+1}\leq \min \left(x_{q},x_{q+1}-1\right),\quad (j<q<l)\]

\item $p(w_{q},x_{q};w_{q+1},x_{q+1})=b_{q}\left(1-b_{q}\right)^{w_{q+1}-w_{q}}$
\item $\ds b_{q}=\frac{\lambda _{q}}{\lambda _{q}+\beta _{q,q+1}}$, $q<N$. 
\end{itemize}
Here $b_{q}$ is the probability of an event that processor~$q$
in state~$x_{q}$ sends at least one message to~$q+1$ before updating
its state $x_{q}\rightarrow x_{q}+1$. For $q=N$ we put $b_{N}=0$.
So in the case $l=N$ the probability~(\ref{eq:pr-x-tran}) takes
the form \[
Z^{-1}\beta _{j,j+1}\, \prod _{q=j+1}^{N-1}p(w_{q},x_{q};w_{q+1},x_{q+1})\, .\]

\paragraph*{Relative coordinates.}

Note that the first processor $x_{1}(t)$ evolves independently of
other processors. It is useful to introduce new process $Y_{c}(t)=(y_{2}(t),\ldots ,y_{N}(t))\in \bZ ^{N-1}$
in relative coordinates as viewing by an observer sitting at the point
$x_{1}(t)$:\[
y_{j}(t):=x_{j}(t)-x_{1}(t),\quad j=2,\ldots ,N\, .\]
In a similar way we define $\Yd{n}=Y_{c}(\tau _{n})$, $n=0,1,\ldots $~.
The free dynamics produce the following transitions of $\Yd{n}$:
\begin{eqnarray}
P\left\{ \Yd{n+1}=y+\be _{j}\, |\, \Yd{n}=y\right\}  & = & \lambda _{j}Z^{-1},\quad j=2,\ldots ,N\, ,\label{eq:trYo}\\
P\left\{ \Yd{n+1}=y-{\textstyle \sum _{j=2}^{N}}\be _{j}\, |\, \Yd{n}=y\right\}  & = & \lambda _{1}Z^{-1}.\label{eq:trYd}
\end{eqnarray}
Since rollback does not affect on the first processor the corresponding
transitions have the same form and the same probabilities as~(\ref{eq:x-tran})
and~(\ref{eq:pr-x-tran}).

\section{Stochastic monotonicity}

All statements of this section are valid for the both Markov processes
$X_{c}^{(N)}(t)$, $t\in \R _{+}$, and $X(n)$, $n\in \bZ _{+}$.
For the sake of breavity we give here results only for the contionuous
time model $X_{c}^{(N)}(t)$. The following results will play a significant
part in the proof of the Theorem~\ref{t:N-3} in Section~\ref{sec:N-3}.

\begin{theorem}\label{t-st-mon}

\noindent Let us consider two cascade models (say $X_{c,1}^{(n)}(t)$
and $X_{c,2}^{(n)}(t)$~) with processors $1,2,\ldots ,n$ and parameters
$\lambda _{1},\ldots ,\lambda _{n}$ and $\beta _{12}^{(1)},\beta _{23}^{(1)},\ldots ,\beta _{n-1,n}^{(1)}$
for the first model $X_{c,1}^{(n)}(t)$ and parameters $\lambda _{1},\ldots ,\lambda _{n}$
and $\beta _{12}^{(2)},\beta _{23}^{(2)},\ldots ,\beta _{n-1,n}^{(2)}$
for the second model $X_{c,2}^{(n)}(t)$. Assume that \[
\beta _{i,i+1}^{(1)}\leq \beta _{i,i+1}^{(2)}\qquad \forall i\, .\]
 Then $X_{c,1}^{(n)}$ is stochastically larger than $X_{c,2}^{(n)}$~,
that is: if $X_{c,1}^{(n)}(0)=X_{c,2}^{(n)}(0)$ then $X_{c,1}^{(n)}(t)\sgeq X_{c,2}^{(n)}(t)$
for any~$t$.~%
\footnote{It means that there exists a \emph{coupling} $\left(\widetilde{X}_{1}^{(n)}(t,\omega ),\widetilde{X}_{2}^{(n)}(t,\omega )\right)$of
stochastic processes $X_{1}^{(n)}(t)$ and $X_{2}^{(n)}(t)$ such
that $P\left\{ \omega :\, \widetilde{X}_{1}^{(n)}(t,\omega )\geq \widetilde{X}_{2}^{(n)}(t,\omega )\, \, \forall t\right\} =1$.
If $w,z\in \R ^{n}$ we say $w\geq z$ if $w_{i}\geq z_{i}$ for all
$i=1,\ldots ,n$~(partial order).%
}

\noindent \end{theorem}

\smallskip{}
\noindent \textbf{Proof} may be given by an explicit coupling construction
of the processes $X_{c,1}^{(n)}(t)$ and $X_{c,2}^{(n)}(t)$ on the
same probability space. The following fact should be used: a~Poisson
flow with intensity $\beta _{12}^{(1)}$ can be obtained from a Poisson
flow with intensity $\beta _{12}^{(2)}$ in which any point (independently
from other) is killed with probability $1-\beta _{12}^{(1)}/\beta _{12}^{(2)}$.
\medskip{}

\begin{cor}[Solid barriers]\label{cor-barrier}

\noindent Fix some $1\leq r_{1}<r_{2}<\cdots <r_{b}<n$ and consider
two cascade models: $X_{c,1}^{(n)}(t)$ with parameters $\left(\, \lambda _{1},\ldots ,\lambda _{n}\, ;\, \beta _{12}^{(1)},\beta _{23}^{(1)},\ldots ,\beta _{n-1,n}^{(1)}\, \right)$
and $X_{c,2}^{(n)}(t)$ with parameters $\left(\, \lambda _{1},\ldots ,\lambda _{n}\, ;\, \beta _{12}^{(2)},\beta _{23}^{(2)},\ldots ,\beta _{n-1,n}^{(2)}\, \right)$,
where\[
\beta _{i,i+1}^{(2)}=\beta _{i,i+1}^{(1)}\quad \forall i\not \in \left\{ r_{1},\ldots ,r_{b}\right\} ,\qquad \beta _{i,i+1}^{(2)}=0\quad \forall i\in \left\{ r_{1},\ldots ,r_{b}\right\} .\]
We can say that the model $X_{c,2}^{(n)}(t)$ differs from the model
$X_{c,1}^{(n)}(t)$ by the presence of $b$ solid barriers between
processors $r_{1}$ and $r_{1}+1$, \ldots, $r_{b}$ and $r_{b}+1$.
Then by Theorem~\ref{t-st-mon} we have that \[
X_{c,1}^{(n)}(t)\sleq X_{c,2}^{(n)}(t)\, .\]
\end{cor}

\section{Case $N=2$}

\label{sec:Case-N-2}

We start with the Markov chain $X_{c}^{(2)}(t)$. Since processor
$1$ works independently, it is enough to consider the Markov chain
$Y_{c}^{(2)}(t)=x_{2}(t)-x_{1}(t)$.

Bearing in mind the remark at the end of Subsection~\ref{sub:Cascade-model},
for brevity of notation let us rescale absolute time in such a way
that $Z=1$. Then the Markov chain $Y(n)$ has the following transition
probabilities\[
p_{i,i+1}=\lambda _{2},\quad p_{i,i-1}=\lambda _{1},\quad p_{i,0}=\beta _{12}\, \, (i\geq 0),\quad p_{i,i}=\beta _{12}\, \, (i<0)\]
and $p_{i,j}=0$ for any another pair $i,j$~.

\begin{center}\includegraphics [width=6in]{n-2.eps}\end{center}

\begin{theorem}\label{t:n-2} If $\lambda _{1}<\lambda _{2}$ then
the Markov chain $\Yd{n}$ is ergodic and we have $\vs _{1}=\vs _{2}=\lambda _{1}$.
If $\lambda _{1}>\lambda _{2}$ then the Markov chain $\Yd{n}$ is
transient and we have $\vs _{1}=\lambda _{1}$, $\vs _{2}=\lambda _{2}$.

\end{theorem}

\begin{proof}The Markov chain $Y(n)$ is one-dimensional and its
analysis is quite easy. To~establish ergodicity under assumption
$\lambda _{1}<\lambda _{2}$ we use the Foster-Lyapunov criterion
(Theorem~\ref{t:Foster}, see~Appendix) with test function $f(y)=|y|$,
$y\in \bZ $. This implies that $x_{2}(t)-x_{1}(t)$ has a limit in
distribution as $t\rightarrow \infty $. Recall that $x_{1}(t)$ is
a Poissonian process hence the limit $\ds \vs _{1}=t^{-1}\lim _{t}x_{1}(t)=\lambda _{1}$
exists (in probability). It follows from this that $\ds \vs _{2}=t^{-1}\lim _{t}x_{2}(t)=\lambda _{1}$.

Under assumption $\lambda _{1}>\lambda _{2}$ we get transience by
choosing the function $f(y)=\min (e^{\delta y},1)$, $y\in \bZ $,
where we fix sufficiently small~$\delta >0$, and applying Theorem~\ref{t:kr-trans}
from Appendix. Therefore any trajectory of $Y(n)$ spends a finite
time in any prefixed domain $\{y\geq C\}$ entailing $\lim _{t\rightarrow \infty }x_{2}(t)-x_{1}(t)=-\infty $
(a.s.). It means that after some time, the messages from $1$ to~$2$
can not produce a rollback anymore, so $x_{1}(t)$ and $x_{2}(t)$
become asymptotically independent and hence $\ds \vs _{2}=t^{-1}\lim _{t}x_{2}(t)=\lambda _{2}$.

\end{proof}

\section{Case $N=3$}

\label{sec:N-3}

\begin{theorem}\label{t:N-3}

Four situations are possible.

\begin{enumerate}
\item If $\lambda _{1}<\min \left(\lambda _{2},\lambda _{3}\right)$ then
$\vs _{1}=\vs _{2}=\vs _{3}=\lambda _{1}$.
\item If $\lambda _{2}>\lambda _{1}>\lambda _{3}$ then $\vs _{1}=\vs _{2}=\lambda _{1}$,
$\vs _{3}=\lambda _{3}$.
\item If $\lambda _{2}<\min \left(\lambda _{1},\lambda _{3}\right)$ then
$\vs _{1}=\lambda _{1}$, $\vs _{2}=\vs _{3}=\lambda _{2}$.
\item If $\lambda _{1}>\lambda _{2}>\lambda _{3}$ then $\vs _{1}=\lambda _{1}$,
$\vs _{2}=\lambda _{2}$, $\vs _{3}=\lambda _{3}$.
\end{enumerate}
\end{theorem}

Items 2, 3 and 4 can be reduced in some sense to the results of the
case $N=2$ (see Theorem~\ref{t:n-2}). We prove them in the current
section. Proof of the item~1 is much more intricate and relies heavily
on the construction of an adequate Lyapunov function needing lengthy
developments deferred to the following section~\ref{sec:N-3-Lyap}.

\paragraph{Proof of Theorem~\ref{t:N-3} (items 2--4).}

We start from the item~2: $\lambda _{2}>\lambda _{1}>\lambda _{3}$.
Since the first two processors are governed by the Markov chain $X_{c}^{(2)}(t)$
and do not depend on the state of processor~3 we apply Theorem~\ref{t:n-2}
and conclude that $X_{c}^{(2)}(t)$ is ergodic and $\vs _{1}=\vs _{2}=\lambda _{1}$.

Let us compare the following two cascade models

\[
X_{c}^{(3)}(t):\qquad 1\stackrel{\beta _{1,2}}{\, \longrightarrow }\, 2\stackrel{\beta _{2,3}}{\, \longrightarrow }\, 3\]

\[
X_{c,2}^{(3)}(t):\qquad 1\stackrel{\beta _{1,2}}{\, \longrightarrow }\, 2\stackrel{0}{\, \longrightarrow }\, 3\]
(parameters $\lambda _{1}$, $\lambda _{2}$ and $\lambda _{3}$ are
the same for the both models $X_{c}^{(3)}(t)$ and $X_{c,2}^{(3)}(t)$~).

In the model $X_{c,2}^{(3)}$ the groups of processors $\{1,2\}$
and $\{3\}$ evolve independently. Evidently, an asymptotic speed
of processor~$3$ in the model $X_{c,2}^{(3)}$ exists and is equal
to~$\lambda _{3}$. By Corollary~\ref{cor-barrier} $X_{c}^{(3)}(t)\sleq X_{c,2}^{(3)}(t)$.
Hence in the model $X_{c}^{(3)}$ an asymptotic speed of the processor~$3$
is \emph{not greater} than $\lambda _{3}$. Since $\lambda _{3}<\lambda _{1}$
we conclude that there exists some time moment $T_{0}$ such that
for $t\geq T_{0}$ in the model~$X_{c}^{(3)}$ messages from $2$
to $3$ that roll back the processor~3 will be very {}``rare''.
So these rare rollbacks will be not essential for an asymptotical
speed of the processor~3. In other words, as $t\rightarrow \infty $
the groups of processors $\{1,2\}$ and $\{3\}$ of the model $X_{c}^{(3)}$
become asymptotically independent, so the processor~3 will move with
the average speed $\lambda _{3}$.

Items 3 and 4 can be considered in a similar way. Note the item~3
consists of two subcases: $\lambda _{1}>\lambda _{3}>\lambda _{2}$
and $\lambda _{3}>\lambda _{1}>\lambda _{2}$. We omit details.

\section{Explicit construction of Lyapunov function}

\label{sec:N-3-Lyap} 

In this section we prove the item~1 of Theorem~\ref{t:N-3}. Recall
that our key assumption here is

\begin{equation}
\lambda _{1}<\lambda _{2},\quad \lambda _{1}<\lambda _{3}.\label{eq:osn-predp}\end{equation}
The main idea is to prove that the Markov chain~$Y(n)$ is ergodic.
To do this we apply the Foster-Lyapunov criterion (see Theorem~\ref{t:Foster}
in Appendix). As in the case of Theorem~\ref{t:n-2} ergodicity of
$Y(n)$ implies that $\vs _{j}=\lambda _{1}$, $j=1,2,3$~.

\subsection{Transition probabilities}

\label{sub:Transition-probabilities}

Consider the embedded Markov chain $Y(n)$. A stochastic dynamics
produced by this Markov chain consists of two components: transitions
generated by the free dynamics and transitions generated by roll-backs.
For each transition probability $p_{\alpha \beta }$, $\alpha \not =\beta $,
we have the following representation:\begin{equation}
p_{\alpha \beta }=s_{\alpha \beta }+r_{\alpha \beta }\, ,\label{eq:p-s-r}\end{equation}
where $s_{\alpha \beta }\geq 0$ corresponds to a transition $\alpha \rightarrow \beta $
which occurs due to the free dynamics and $r_{\alpha \beta }$ corresponds
to a roll-back transition $\alpha \rightarrow \beta $. 

Taking into account the remark at the end of Subsection~\ref{sub:Cascade-model},
without loss of generality we assume that the time is rescaled in
such way that $Z=1$. This slightly simplifies notation for transition
probabilities. For example, free dynamics transitions~(\ref{eq:trYo})--(\ref{eq:trYd})
are equal to $\lambda _{2}$, $\lambda _{3}$ and $\lambda _{1}$
correspondingly. On the next figure we show all non-zero transitions
$\alpha \rightarrow \beta $, $(\alpha \not =\beta )$. It is true,
of course, that $p_{\alpha \alpha }=1-\sum _{\beta \not =\alpha }p_{\alpha \beta }$,
but it is useless to put this information on the picture. Below we
give the explicit form of rollback transition probabilities: \begin{eqnarray*}
1\rightarrow 2: & \qquad  & r_{yz}=\beta _{12}\quad \quad \textrm{for }0<y_{2}\\
2\rightarrow 3: &  & r_{yz}=\beta _{23}\quad \quad \textrm{for }y_{2}<y_{3}\\
1\rightarrow 2\rightarrow 3: &  & r_{yz}=\left\{ \begin{array}{rl}
 \beta _{12}\left(1-b_{2}\right)^{z_{3}}b_{2},\quad  & z_{3}<y_{3}\\
 \beta _{12}\left(1-b_{2}\right)^{y_{3}},\quad  & z_{3}=y_{3}\end{array}\right.\quad \textrm{for }0<y_{3}\leq y_{2}\\
1\rightarrow 2\rightarrow 3: &  & r_{yz}=\left\{ \begin{array}{rl}
 \beta _{12}\left(1-b_{2}\right)^{z_{3}}b_{2},\quad  & z_{3}\leq y_{2}\\
 \beta _{12}\left(1-b_{2}\right)^{y_{2}+1},\quad  & z_{3}=y_{3}\end{array}\right.\quad \textrm{for }0<y_{2}<y_{3}
\end{eqnarray*}
 where $y=(y_{2},y_{3})$, $z=(z_{2},z_{3})$.

\includegraphics [width=5in]{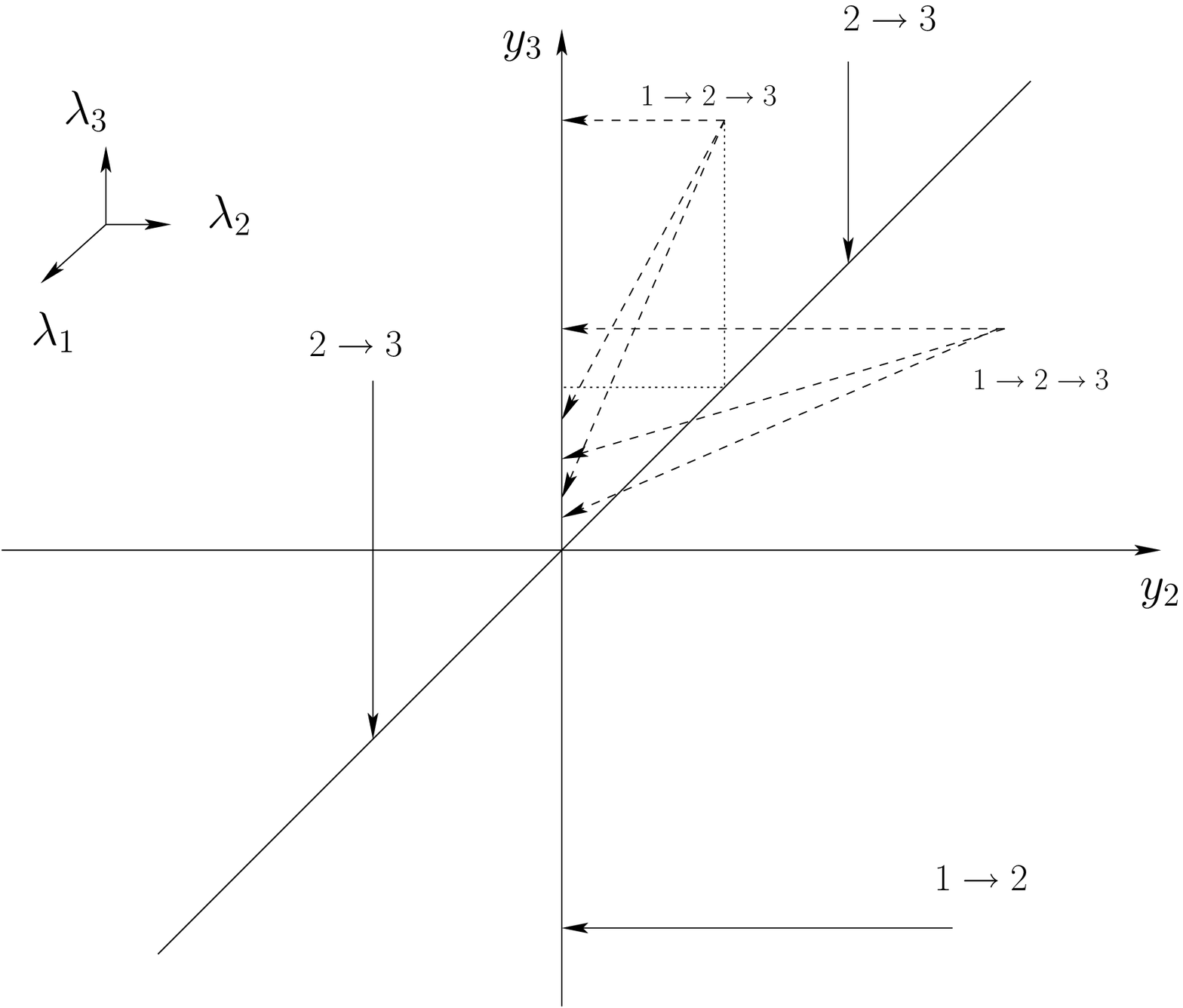}

\subsection{Contour of level 1}

In the plane $Oy_{2}y_{3}$ consider the ellipse\[
e(y_{2},y_{3})=ay_{2}^{2}+b\left(y_{2}-y_{3}\right)^{2}=1,\qquad a>0,b>0,\]
and draw a tangent line to it with normal vector $\left(-\Delta ,1\right)$.
Evidently, there exist two tangent lines with the same normal vector
$\left(-\Delta ,1\right)$. If $\Delta >0$ is sufficiently large
then one of this tangent line touches the ellipse at some point $T_{3}$
of the domain $y_{2}<0$, $y_{3}<0$. Take a segment on this line
from the point $T_{3}$ to a point $K_{3}=(0,u_{3})$ of intersection
with coordinate axis $Oy_{3}$. Now let us draw tangent lines to the
ellipse corresponding to a normal vector $\left(1,-\Delta \right)$.
If $\Delta >0$ is sufficiently large, then one of these lines touches
the ellipse at some point $T_{2}$ of the domain $y_{3}<0$. Let us
take this tangent line and fix a segment on it from the point $T_{2}$
to a point $K_{2}=(u_{2},0)$ of intersection with coordinate axis~$Oy_{2}$.
It is evident that $[K_{2}K_{3}]=\R _{+}^{2}\cap \left\{ (y_{2},y_{3}):\, y_{2}/u_{2}+y_{3}/u_{3}=1\right\} $. 

Let us consider now a closed contour $L$, consisting of subsequently
joined segment $K_{3}K_{2}$, segment $K_{2}T_{2}$, arc $T_{2}T_{3}$
of the ellipse and segment $T_{3}K_{3}$. This contour has the following
property: any ray of the form$\{cv,\, c>0\}$, where $v\in \R ^{2}$,
$v\not =0$, has exactly one common point with the contour $L$. 

\begin{center}\includegraphics [width=5in]{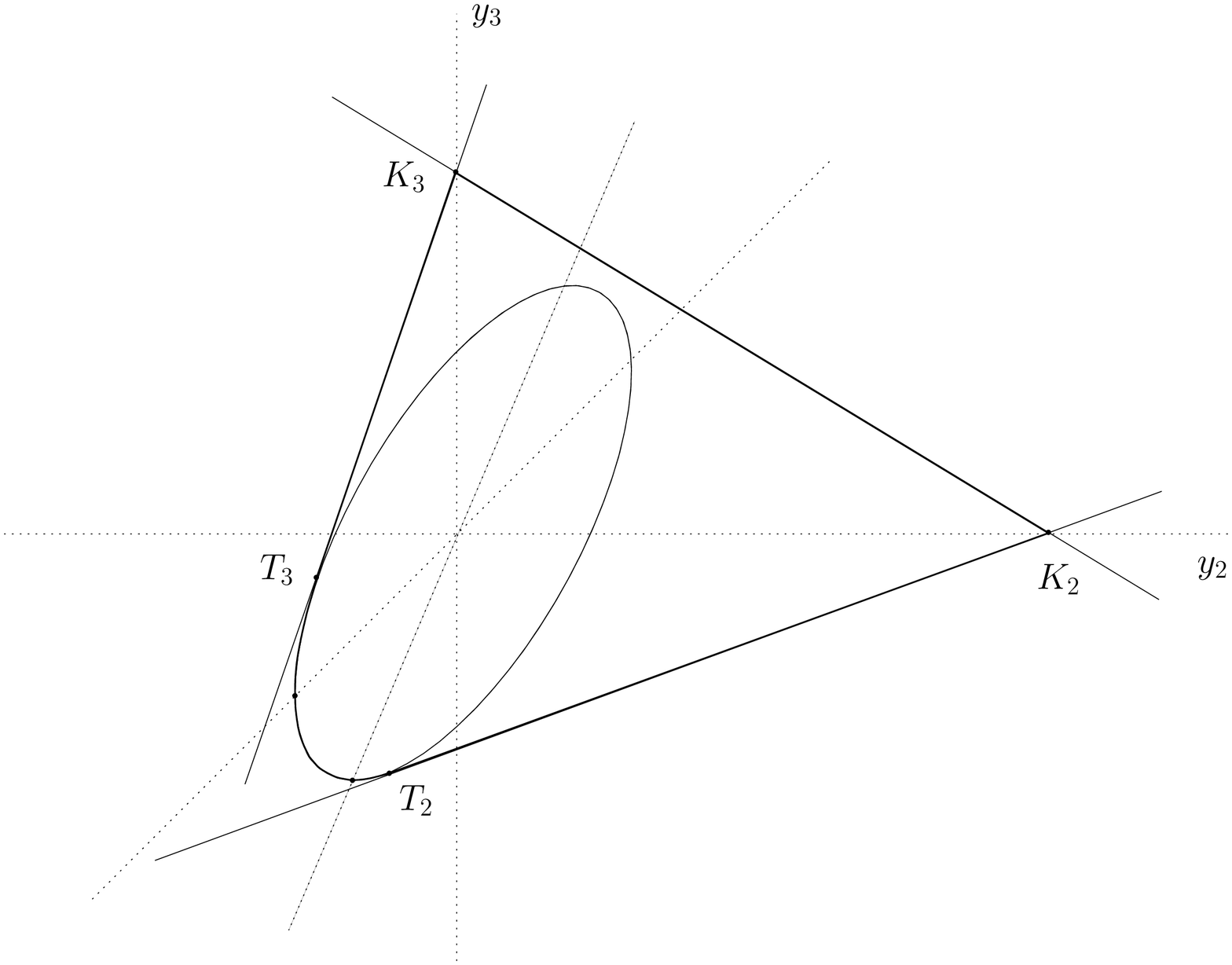}\end{center}

We denote by $n(y)$ the outer normal unitary vector of the contour
$L$ corresponding to the point $y\in L$, $n(y)$ is well defined
at all points of~$L$ except the points $K_{2}$ and $K_{3}$ and,
moreover, this function is continuous on~$L$ except the points~$K_{2}$
and $K_{3}$. The~behaviour of $n(y)$ on the arc $T_{2}T_{3}$ is
of prime interest:

\[
n(y)=\frac{\nabla e(y)}{\left\Vert \nabla e(y)\right\Vert },\qquad \nabla e(y)=2(\, ay_{2}+b(y_{2}-y_{3}),-b(y_{2}-y_{3})\, ),\quad y\in T_{2}T_{3}\subset L\, .\]
It is easy to see that $n(y)=n(T_{2})$ for $y\in (K_{2}T_{2}]$,
$n(y)=n(T_{3})$ for $y\in [T_{3}K_{3})$ and\[
n(y)=\left(u_{2}^{-1},u_{3}^{-1}\right)\qquad y\in (K_{3}K_{2}).\]

For the sequel it is important to point out the following points of
the arc $T_{2}T_{3}$: $y^{(3)}=(-a^{-1/2},-a^{-1/2})$ and $y^{(2)}$,
$\left\{ y^{(2)}\right\} =T_{2}T_{3}\cap \left\{ y_{3}^{(2)}=\frac{a+b}{b}y_{2}^{(2)}\right\} $.
It is easy to check that \[
n(y^{(2)})\, \Vert \, Oy_{3},\quad n(y^{(3)})\, \Vert \, Oy_{2}\, .\]
Obviously, both points belong to the domain $\left\{ y_{2}<0,\, y_{3}<0\right\} $.

\begin{lemma} ~\label{l-normal}

The function $n(y)$ has the following properties:

\begin{itemize}
\item $\skp{n(y)}{y}\not =0$ $\forall y\in L\backslash \{K_{2},K_{3}\}$
\item If $\Delta >0$ is sufficiently large then there exist continuous
functions $c_{2}(y)$ and $c_{3}(y)$ such that \[
c_{2}(T_{2})=c_{3}(T_{3})=1,\qquad c_{2}(T_{3})=c_{3}(T_{2})=0,\qquad c_{2}(y)>0,\, c_{3}(y)>0\quad y\in (T_{2},T_{3})\]
\[
n(y)=c_{2}(y)n(T_{2})+c_{3}(y)n(T_{3}),\qquad y\in (T_{2},T_{3}).\]

\item $\skp{n(y)}{(0,-1)}<0$ if $y_{2}<0$, $y_{3}>y_{2}$, and $\skp{n(y)}{(-1,0)}<0$
if $y_{2}>0$, $y_{3}<0$.
\end{itemize}
\end{lemma}

\subsection{Definition of function $\flam $}

For any point $(y_{2},y_{3})\in \R ^{2}\backslash \{0\}$ define $\flam (y_{2},y_{3})>0$
such that \[
\frac{(y_{2},y_{3})}{\flam (y_{2},y_{3})}\, \in \, L\, .\]
For $(y_{2},y_{3})=0$ we put $\flam (0,0)=0$. The function $\flam (y_{2},y_{3})$
is well-defined and has the following properties:

\begin{itemize}
\item $\flam :\, \R ^{2}\rightarrow \R _{+}$ (positivity)
\item $\flam (ry_{2},ry_{3})=r\flam (y_{2},y_{3})$, $r>0$, (homogeneity)
\item $L=\{y:\, \flam (y)=1\}$.
\end{itemize}
To any point $y=(y_{2},y_{3})$ we put in correspondance a point $y^{*}:=\frac{y}{\flam (y)}\in L$.
Therefore, $\flam (y^{*})=1$.

\begin{center}\includegraphics [width=5in]{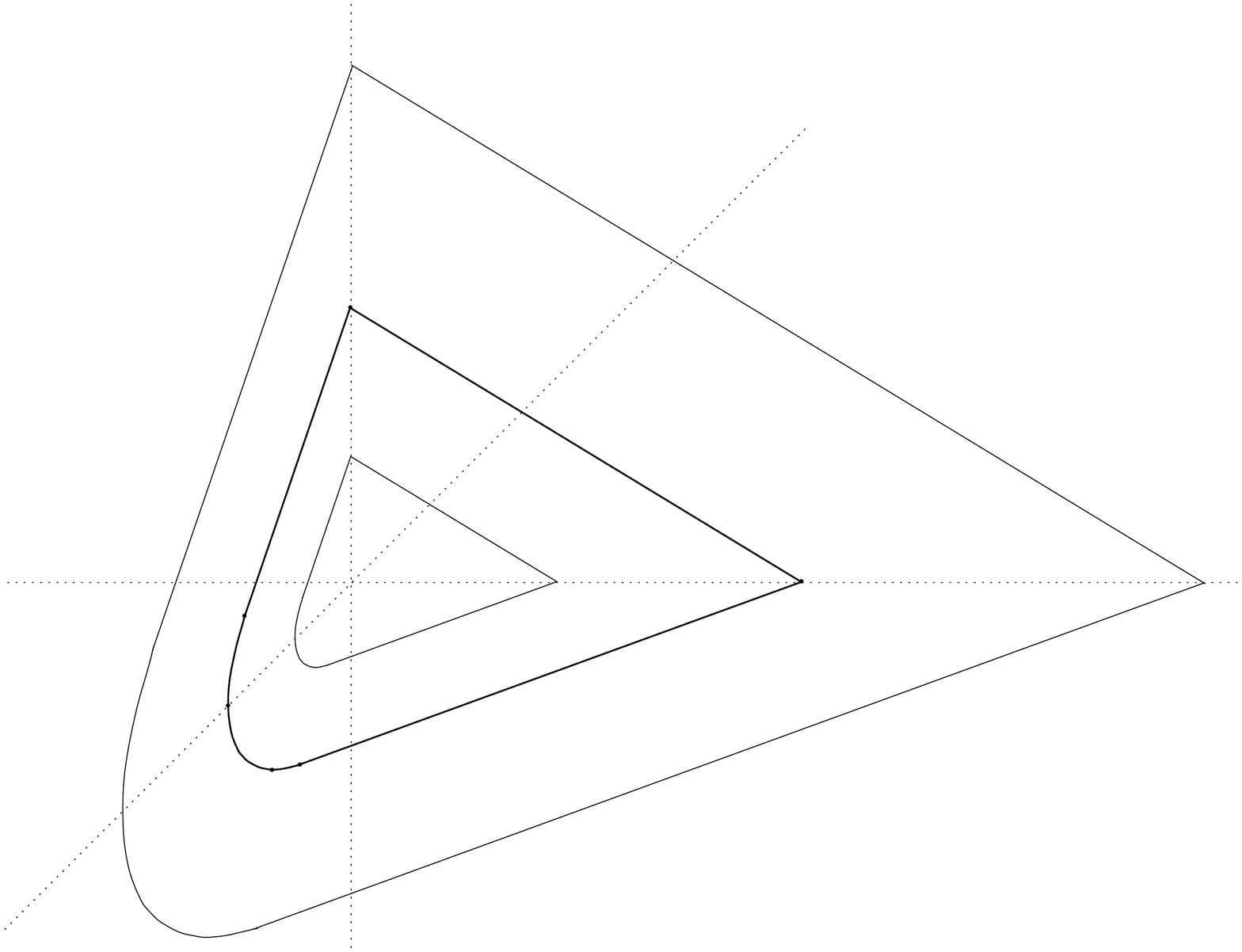}\end{center}

\begin{lemma}

\label{l-grad-fi}~

\begin{itemize}
\item The gradient $\nabla \flam (y)$ exists at all points except that
$y$ for which $y^{*}=K_{2}$ or $K_{3}$ and, moreover, the gradient
is constant on rays of the form $\{cv,\, c>0\}$, $v\in \R ^{2}$:\[
\nabla \flam (y)=\frac{n(y^{*})}{\skp{y^{*}}{n(y^{*})}}\, .\]

\item Let $y=(y_{2},y_{3})$ be such that $y^{*}\in T_{2}T_{3}$. Then\begin{equation}
\left|\flam (w)-\skp{\nabla \flam (y^{*})}{w}\right|\leq \frac{const}{\flam (y)}\, \left\Vert w-y\right\Vert ^{2}.\label{eq:lagrange-ineq}\end{equation}
In other words, in a neighbourhood of the point $y$ the function
$\flam $ can be approximated by the linear function $\skp{\nabla \flam (y^{*})}{\cdot }$.
\end{itemize}
\end{lemma}

In particular, $\flam (y)=\skp{\nabla \flam (y^{*})}{y}$.

Proof of Lemma~\ref{l-grad-fi} is a straightforward computation.

\subsection{Modification of the principle of local linearity}

For any state $\alpha $ define a set $T_{\alpha }=\{\beta :\, p_{\alpha \beta }>0\}$.
Recall decomposition~(\ref{eq:p-s-r}) and define $F_{\alpha }=\{\beta :\, s_{\alpha \beta }>0\}$
and $R_{\alpha }=\{\beta :\, r_{\alpha \beta }>0\}$. It is evident
that $T_{\alpha }=F_{\alpha }\cup R_{\alpha }$. The most simple case
is $F_{\alpha }\cup R_{\alpha }=\varnothing $. The case $F_{\alpha }\cup R_{\alpha }\not =\varnothing $
can be reduced to the previous one by a dilatation of the state space.
Thus we assume that $F_{\alpha }\cup R_{\alpha }=\varnothing $ and
consider the events $\left\{ \Xc (n+1)\in F_{\alpha }\right\} $ and
$\left\{ \Xc (n+1)\in R_{\alpha }\right\} $. On the set $\left\{ \omega \in \Omega :\, \Xc (n,\omega )=\alpha \right\} $
we have $\IF{\alpha }(\omega )+\IR{\alpha }(\omega )\equiv 1$. Hence,
\begin{eqnarray*}
\sE \left(f(\Xc (n+1))\, |\, \Xc (n)=\alpha \right)-f(\alpha ) & = & \sE \left(\left(f(\Xc (n+1))-f(y)\right)\IF{\alpha }\, |\, \Xc (n)=\alpha \right)+\\
 &  & +\sE \left(\left(f(\Xc (n+1))-f(y)\right)\IR{\alpha }\, |\, \Xc (n)=\alpha \right)
\end{eqnarray*}

It follows from definition of the Markov chain~$Y(n)$ (see Subsect.~\ref{sub:Transition-probabilities})
that the diameters $d_{\alpha }:=diam\, F_{\alpha }$ are uniformly
bounded in~$\alpha $: $d=\max _{\alpha }d_{\alpha }<+\infty $.
Define a vector\[
M_{F}(\alpha )=\sE \left(\left(\Xc (n+1)-\alpha \right)\IF{\alpha }\, |\, \Xc (n)=\alpha \right)=\sum _{\beta \in F_{\alpha }}(\beta -\alpha )p_{\alpha \beta }\, .\]
This is an analogue of a notion of mean jump (see~(\ref{eq:M-alpha})
in Appendix). In the next subsection we shall need the following modification
of the principle of local linearity from~\cite{FMM} (see also Subsect.~\ref{sub:Principle-of-local}
in Appendix).

\begin{lemma}\label{l-mod-loclin}Assume that the following condition
holds \[
\inf _{l}\, \sup _{\tilde{\alpha }\in \R ^{n},\| \tilde{\alpha }-\alpha \Vert \leq d_{\alpha }}\, \left|f(\tilde{\alpha })-l(\tilde{\alpha })\right|\, <\, \varepsilon \, ,\]
where $\inf $ is taken over all linear functions~$l$. If \[
f\left(\alpha +M_{F}(\alpha )\right)-f(\alpha )<-5\varepsilon \, ,\]
then the following inequality\[
\sE \left(\, \left(\, f(\Xc (n+1))-f(\alpha )\, \right)\IF{\alpha }\, |\, \Xc (n)=\alpha \, \right)<-\varepsilon \, \]
 holds.

\end{lemma}

The proof of this statement repeats the proof of principle of local
linearity presented in~\cite{FMM} and is omitted.

\subsection{Proof of the Foster condition}

The validity of the Foster condition will follow from several ancillary
lemmas dealing with the following different domains of the state space:
\begin{eqnarray*}
E_{-} & := & \left\{ y=(y_{2},y_{3}):\, \min (y_{2},y_{3})<0\right\} ,\\
E_{1} & := & \left\{ y=(y_{2},y_{3}):\, y_{2}>0,y_{3}>0\right\} ,\\
E_{1,2} & := & \{y_{2}>0,\, y_{3}=0\},\\
E_{1,3} & := & \{y_{2}=0,\, y_{3}>0\}.
\end{eqnarray*}

\begin{lemma}\label{l-leftlower} Consider the domain $E_{-}=\left\{ y=(y_{2},y_{3}):\, \min (y_{2},y_{3})<0\right\} $.
There exists $\CC{l-leftlower}>0$, such that if $\flam (y)>\CC{l-leftlower}$,
then

~\quad ~a) ~\[
\sE \left(\left(\flam (\Xc (n+1))-\flam (y)\right)\IR{y}\, |\, \Xc (n)=y\right)\leq 0\, ,\]

~\quad ~b) there exists $\varepsilon >0$ such that\[
\flam \left(y+M_{F}(y)\right)-\flam (y)<-5\varepsilon \, .\]

\end{lemma}

\begin{proof}

It is evident that the vector\[
M_{F}(y)=\left(\lambda _{2}-\lambda _{1},\lambda _{3}-\lambda _{1}\right)\]
is constant (does not depend on $y$). Since the vector $n(T_{2})$
is co-directed with the vector $\left(1,-\Delta \right)$ and the
vector $n(T_{3})$ is co-directed with the vector $\left(-\Delta ,1\right)$
and the conditions $\lambda _{2}>\lambda _{1}$, $\lambda _{3}>\lambda _{1}$
hold, we can find a large $\Delta _{1}>0$ such that \[
\skp{M_{F}(y)}{n(T_{2})}<0,\quad \skp{M_{F}(y)}{n(T_{3})}<0\, ,\qquad \forall \Delta >\Delta _{1}\, .\]
Fix this $\Delta _{1}$. Hence, by Lemma~\ref{l-normal} there exists
$\varepsilon >0$ such that\[
\skp{M_{F}(y)}{n(y)}<-6\varepsilon \quad \textrm{if }\, \min (y_{2},y_{3})<0.\]
Put $w=y+M_{F}(y)$ and consider\begin{eqnarray*}
\flam (w)-\flam (y) & = & \flam (w)-\skp{\nabla \flam (y^{*})}{y}\\
 & = & \flam (w)-\skp{\nabla \flam (y^{*})}{w}+\skp{\nabla \flam (y^{*})}{w-y}\\
 & = & \flam (w)-\skp{\nabla \flam (y^{*})}{w}+\skp{\nabla \flam (y^{*})}{M_{F}(y)}\, .
\end{eqnarray*}
By~(\ref{eq:lagrange-ineq}) for any given $\varepsilon >0$ we can
choose $C_{0}>0$ such that \[
\left|\flam (y+M_{F}(y))-\skp{\nabla \flam (y^{*})}{M_{F}(y)}\right|\leq \frac{const}{\flam (y)}\, d^{2}\leq \varepsilon \, \qquad \textrm{if }\flam (y)\geq C_{0}\, .\]
Now the item b) of the lemma easily follows.

Let us prove the item a) of the lemma. Note that in the domain $y_{3}>y_{2},\, y_{2}<0$
a rollback decreases coordinate $y_{3}$: $(y_{2},y_{3})\rightarrow \left(y'_{2},y'_{3}\right)=(y_{2},y_{2})$.
From geometrical properties of level sets of function $\flam $ and
item~3 of Lemma~\ref{l-normal} it follows that any transition generated
by a rollback decreases a value of the function~$\flam $: $\flam (\, (y_{2},y_{3})\, )<\flam (\, (y'_{2},y'_{3})\, )$.
In the domain $y_{3}<0,\, y_{2}>0$ a rollback has the following form:
$(y_{2},y_{3})\rightarrow \left(y'_{2},y'_{3}\right)=(0,y_{3})$.
For similar reasons we again have $\flam (\, (y_{2},y_{3})\, )<\flam (\, (y'_{2},y'_{3})\, )$.
In the domain $y_{3}\leq y_{2}<0$ there is no rollback. Now the item
a) easily follows.\end{proof}

\begin{lemma}\label{l-posoct}

Consider the domain: $E_{1}=\left\{ y=(y_{2},y_{3}):\, y_{2}>0,y_{3}>0\right\} $. 

\begin{enumerate}
\item The conditional expectation\[
\sE \left(\left(\flam (\Xc (n+1))-\flam (y)\right)\IF{y}\, |\, \Xc (n)=y\right)=\skp{\left(u_{2}^{-1},u_{3}^{-1}\right)}{M_{F}(y)}\]
does not depend on $y$.
\item There exist constants $\CC{l-posoct},\ga{l-posoct}>0$ such that \begin{equation}
\sE \left(\left(\flam (\Xc (n+1))-\flam (y)\right)\IR{y}\, |\, \Xc (n)=y\right)\leq -\ga{l-posoct}\flam (y)\, \quad \textrm{if}\, \textrm{ }\flam (y)>\CC{l-posoct}\label{eq:pos-quat-ineq}\end{equation}

\end{enumerate}
\end{lemma}

\begin{proof}

The first statement follows from the fact that in this domain $\flam (y)=\skp{\left(u_{2}^{-1},u_{3}^{-1}\right)}{y},$
and the vector $M_{F}(y)$ does not depend on $y$. 

Let us prove~(\ref{eq:pos-quat-ineq}). Fix some level set \[
L_{C}^{+}:=\left\{ y:\, \flam (y)=C\right\} \cap E_{1}\equiv \left\{ y_{2}/u_{2}+y_{3}/u_{3}=C,\, y_{2}>0,y_{3}>0\right\} \]
 and consider an action of rollbacks for $y\in L_{C}^{+}$. We have
three different situations.

a) Let $y$ be such that $y_{2}\geq y_{3}>0$. It follows that $\ds y_{2}\geq \frac{C}{\frac{1}{u_{2}}+\frac{1}{u_{3}}}\, .$
As it can be easily concluded from Subsection~\ref{sub:Transition-probabilities},
with probability $\beta _{12}$ we have a rollback of the following
form $(y_{2},y_{3})\rightarrow (0,y_{3}')$ where $0\leq y_{3}'\leq y_{3}$.
Then we obtain\begin{eqnarray*}
\flam \left((0,y_{3}')\right)-\flam \left((y_{2},y_{3})\right) & = & \left(\frac{0}{u_{2}}+\frac{y_{3}'}{u_{3}}\right)-\left(\frac{y_{2}}{u_{2}}+\frac{y_{3}}{u_{3}}\right)\\
 & \leq  & -\frac{y_{2}}{u_{2}}\, \leq \, -\frac{C}{\frac{1}{u_{2}}+\frac{1}{u_{3}}}\, \, \left(u_{2}\right)^{-1}=-\frac{C}{1+u_{2}/u_{3}}\, 
\end{eqnarray*}
uniformly in $y_{3}'$ such that $y_{3}'\leq y_{3}$. To phrase it,
we will say that with probability~$\beta _{12}$ the increment of
$\flam (y)$ is less or equal to $\, -\, \ds \frac{C}{1+u_{2}/u_{3}}$~.
Hence the conditional mean \begin{equation}
\sE \left(\left(\flam (\Xc (n+1))-\flam (y)\right)\IR{y}\, |\, \Xc (n)=y\right)\label{eq:c-mean}\end{equation}
 does not exceed the value $\, -\ds \frac{\beta _{12}C}{1+u_{2}/u_{3}}$
if $y\in L_{C}^{+}$, $y_{2}\geq y_{3}>0$.

b) Let $y\in L_{C}^{+}$ be such that $0<\frac{1}{2}y_{3}\leq y_{2}<y_{3}$.
It follows that $y_{2}\geq \frac{C}{2\left(\frac{1}{u_{2}}+\frac{1}{u_{3}}\right)}$.
With probability $\beta _{12}$ we have a rollback $(y_{2},y_{3})\rightarrow (0,y_{3}')$
where $0\leq y_{3}'\leq y_{3}$ and with probability~$\beta _{23}$
we have a rollback $(y_{2},y_{3})\rightarrow (y_{2},y_{2})$. Both
of them give negative increments of the function $\flam $. But the
first rollback gives the increment $\flam \left((0,y_{3}')\right)-\flam \left((y_{2},y_{3})\right)$
which is less or equal to~$\, -\, \ds \frac{C}{2(1+u_{2}/u_{3})}$.
So we conclude that the above conditional mean~(\ref{eq:c-mean})
will not exceed the value $\, -\frac{1}{2}\beta _{12}\left(1+u_{2}/u_{3}\right)^{-1}C$.

c) Now let $y\in L_{C}^{+}$ be such that $0<y_{2}\leq \frac{1}{2}y_{3}$.
It follows that $y_{3}-y_{2}\geq K(C)$ where\[
K(C):=\frac{1}{2}\cdot \frac{C}{\left(\frac{1}{2u_{2}}+\frac{1}{u_{3}}\right)}=\frac{Cu_{3}}{u_{3}/u_{2}+2}\, .\]
 With probability $\beta _{12}$ we have a rollback $(y_{2},y_{3})\rightarrow (0,y_{3}')$,
$0\leq y_{3}'\leq y_{3}$, and with probability~$\beta _{23}$ we
have a rollback $(y_{2},y_{3})\rightarrow (y_{2},y_{2})$. The first
rollback gives a negative increment of the function~$\flam $, and
the second rollback gives the increment $\flam \left((y_{2},y_{2})\right)-\flam \left((y_{2},y_{3})\right)$
which is less or equal to $\, -\, K(C)/u_{3}$. Hence the conditional
expectation~(\ref{eq:c-mean}) does not exceed the value $\, -\beta _{23}K(C)/u_{3}=-\beta _{23}\left(u_{3}/u_{2}+2\right)^{-1}C$.

The proof of the lemma is completed.

\end{proof}

\begin{lemma}\label{l-axes}

Consider the cases when $y$ belongs to the axes: $y\in E_{1,3}$,
$y\in E_{2,3}$. Here\[
\sE \left(\left(\flam (\Xc (n+1))-\flam (y)\right)\IR{y}\, |\, \Xc (n)=y\right)\leq -\ga{l-axes}\flam (y)\, ,\]
and\[
\sE \left(\left(\flam (\Xc (n+1))-\flam (y)\right)\IF{y}\, |\, \Xc (n)=y\right)\]
does not depend on $y$, where $y\in \{y:\, \flam (y)>\CC{l-axes}\}$.

\end{lemma}

\begin{proof}

We consider in details the case $E_{1,3}=\{y_{2}=0,\, y_{3}>0\}$.
We start with the free dynamics. The following transition \[
(0,y_{3})\rightarrow (y_{2}',y_{3}')=(-1,y_{3}-1)\in E_{3}=\{y_{2}<0,y_{3}>y_{2}\}.\]
occurs with probability $\nconst \lambda _{1}$. It is easy to see
that for $y_{3}>\CC{l-axes}$ values of the function $\flam (\cdot )$
in both points $(0,y_{3})$ and $(-1,y_{3}-1)$ concide with the values
of linear function $\skp{n(T_{3})}{\cdot }$. 

With probability $\nconst \lambda _{2}$ we have a transition\[
(0,y_{3})\rightarrow (y_{2}',y_{3}')=(1,y_{3})\in E_{1}=\{y_{2}>0,y_{3}>0\},\]
and with probability $\nconst \lambda _{3}$ we have a transition\[
(0,y_{3})\rightarrow (y_{2}',y_{3}')=(0,y_{3}+1)\in E_{3}.\]
Evidently, that in $(0,y_{3})$, $(1,y_{3})$ and $(0,y_{3}+1)$ the
values of $\flam (\cdot )$ coincide with the values of a linear function
$\skp{n_{1}}{\cdot }$. Hence,

\begin{eqnarray*}
\sE \left(\left(\flam (\Xc (n+1))-\flam (y)\right)\IF{y}\, |\, \Xc (n)=y\right) & = & \lambda _{1}\skp{n(T_{3})}{(-1,-1)}\\
 &  & \null +\lambda _{2}\skp{n_{1}}{(1,0)}+\lambda _{3}\skp{n_{1}}{(0,1)}.\\
 &  & 
\end{eqnarray*}
Since the r.h.s.\ does not depend on~$y$ we get the second statement
of the lemma. 

Due to a rollback the Markov chain $Y(n)$ goes from the state $(0,y_{3})\in E_{1,3}$
to a state $(0,0)$ with probability $\nconst \beta _{23}$. Note
that values of $\flam (\cdot )$ at these two points can be calculated
by using the linear function $\skp{n_{1}}{\cdot }$. Obviously, that
the increment of~$\flam $ corresponding to this rollback is equal
to $\flam \left(\, (0,0)\, \right)-\flam \left(\, (0,y_{3})\, \right)=-C$,
where $C=\flam (\, (0,y_{3})\, )$. The first statement of the lemma
is proved.

The case of the domain $E_{1,3}=\{y_{2}>0,\, y_{3}=0\}$ is similar.

\end{proof}

\begin{lemma}\label{l-origin} For any $C_{0}>0$\[
\sup _{y:\, \flam (y)\leq C_{0}}\sE \left(\flam (\Xc (n+1))\, |\, \Xc (n)=y\right)<+\infty \, .\]

\end{lemma}

\begin{proof}

This statement follows from the fact that the jumps of any fixed neighbourhood
of $(0,0)$ are bounded and the fact that the function $\flam $ is
continuous. \end{proof}

In view of the Lemmas \ref{l-mod-loclin}--\ref{l-origin}, the Foster-Lyapunov
criterion (Theorem~\ref{t:Foster}, Appendix) is fulfilled with $f(y)=\flam (y)$
therefore the Markov chain~$Y(n)$ is ergodic and hence the proof
of the item~1 of Theorem~\ref{t:N-3} is completed.

\section{Conclusions, conjectures and perspectives}

\label{sec:Conclusions}

\subsection{Decomposition into groups}

We shall always assume that all $\lambda _{1},\ldots ,\lambda _{N}$
are different. Define a function \[
\ell (m):=\min _{i\leq m}\lambda _{i}\, .\]
Evidently, this function has the following property:\[
\ell (1)=\lambda _{1}\geq \cdots \geq \ell (m)\geq \ell (m+1)\geq \cdots \geq \min (\lambda _{1},\ldots ,\lambda _{N})=\ell (N)\, .\]

Level sets of function~$\ell $ generate a partition of the set $\left\{ 1,\ldots ,N\right\} $.
Namely, there exists a sequence $j_{1}=1<j_{2}<\cdots <j_{K}<j_{K+1}=N+1$
such that the set of all processors can be divided into several nonintersecting
groups \begin{equation}
\left\{ 1,\ldots ,N\right\} =\bigcup _{k=\overline{1,K}}G_{k}\, ,\label{eq:part-Gk}\end{equation}
 \[
G_{k}:=\, \left\{ j_{k},j_{k}+1,\ldots ,j_{k+1}-1\right\} ,\, \quad \, \ell (j_{k}-1)>\ell (j_{k})=\cdots =\ell (j_{k+1}-1)>\lambda _{j_{k+1}}\, .\]

\begin{rem}

\noindent \emph{An equivalent description of the group is possible.
We say, for example, that $1,2,\ldots ,k$ is a group if \begin{equation}
\lambda _{1}\leq \min \left(\lambda _{2},\ldots ,\lambda _{k}\right),\quad \lambda _{1}>\lambda _{k+1}.\label{eq:char-group}\end{equation}
}

\end{rem}

\subsection{Long-time behaviour of the groups}

Taking into account Theorems~\ref{t:n-2} and~\ref{t:N-3} and the
above notion of groups of processors we put forward the following
Conjecture.

\noindent \begin{conj}\label{con-N-main} \emph{}Assume that all
$\lambda _{1},\ldots ,\lambda _{N}$ are different. For any $j$ the
following limit $\ds \vs _{j}=\lim _{t\rightarrow +\infty }\frac{x_{j}(t)}{t}$
exists and $\vs _{j}=\ell (j)$ .

\noindent \end{conj}

Therefore this conjecture entails $\vs _{j}=\ell (j_{k})$ for $j\in G_{k}$.
If for some $k$ the group $G_{k}$ consists of \emph{more than one}
processor we may say that the processors of the group $G_{k}$ are
\emph{synchronized}.

\begin{rem}[On monotone cases]

\noindent \emph{If $\lambda _{1}<\cdots <\lambda _{N}$ then $\vs _{j}=\lambda _{1}$
for any~ $j$. \\ If $\lambda _{1}>\cdots >\lambda _{N}$ then for
all~ $j$ we have $\vs _{j}=\lambda _{j}$. }
\smallskip{}

\end{rem}

Let us discuss briefly perspectives of rigorous proof of the above
Conjecture for large values of~$N$. In fact, we have already proved
this conjectures for a wide class of cascade models.

\begin{theorem} 

Assume that all $\lambda _{1},\ldots ,\lambda _{N}$ are different
and a partition~(\ref{eq:part-Gk}) of the set of processors $\left\{ 1,\ldots ,N\right\} $
is such that $\left|G_{k}\right|\leq 3$ for all~$k$. Then the limits
$\ds \vs _{j}=\lim _{t\rightarrow +\infty }\frac{x_{j}(t)}{t}$ exist
and $\vs _{j}=\ell (j)$ .

\end{theorem}

The proof of this statement is just a combination of the result of
Theorem~\ref{t:N-3} (item~1) and arguments of the proof of items
2-4 of Theorem~\ref{t:N-3}. We will not pursue further.

So the key to the proof of Conjecture consists in generalization of
item~1 of Theorem~\ref{t:N-3}. As it was seen in Section~\ref{sec:N-3-Lyap},
a possible way of such generalization is an explicit construction
of Foster-Lyapunov function in high dimensions. This seems to be a
difficult technical problem which is out of scope of this paper.

\appendix

\section{Appendix}

Let $\left(\xi _{n},\, n=0,1,\ldots \right)$ be a countable irreducible
aperiodic Markov chain with the state space $\mathcal{A}$.

\subsection{Criteria}

We use the following Foster criterion.

\begin{theorem}[\cite{FMM}]\label{t:Foster}The Markov chain $\xi _{n}$
is ergodic if and only if there exists a positive function $f(\alpha )$,
$\alpha \in \mathcal{A}$, a number $\varepsilon >0$ and a finite
set $A\in \mathcal{A}$ such that

1) \[
\sE \left(f(\xi _{n+1})\, |\, \xi _{n}=y\right)-f(y)<-\varepsilon \]
for all $y\not \in A$,

2) $\sE \left(f(\xi _{n+1})\, |\, \xi _{n}=y\right)<+\infty $ for
all $y\in A$.

\end{theorem}

The following theorem give a criterion of transience.

\begin{theorem}[\cite{FMM}]\label{t:kr-trans} The Markov chain is
transient, if and only if there exists a positive function $f(\alpha )$
and a set $A$ such that the following inequalities are fulfilled\[
\sE \left(f(\xi _{m+1})\, |\, \xi _{m}=\alpha _{i}\right)-f(\alpha _{i})\leq 0,\quad \forall \alpha _{i}\not \in A,\]
\[
f(\alpha _{k})<\inf _{\alpha _{j}\in A}f(\alpha _{j}),\quad \textrm{for at least one }\alpha _{k}\not \in A.\]

\end{theorem}

\subsection{Principle of local linearity}

\label{sub:Principle-of-local}

From now on we assume that the state space $\mathcal{A}$ is some
subset of $\R ^{k}$. Define a vector of mean jump from the point
$\alpha $\begin{equation}
M(\alpha )=\sE \left(\, \left(\xi _{n+1}-\alpha \right)\, |\, \xi _{n}=\alpha \right)=\sum _{\beta }(\beta -\alpha )p_{\alpha \beta }\, .\label{eq:M-alpha}\end{equation}
Assume that $d_{\alpha }:=\max _{\beta }\left\{ \left|\beta -\alpha \right|:\, p_{\alpha \beta }>0\right\} <\infty $
for all $\alpha $.

The following principle of local linearity was proved in~\cite{FMM}. 

\begin{lemma}Assume that at some point $\alpha $ the following condition\[
\inf _{l}\, \sup _{\tilde{\alpha }\in \R ^{n},\| \tilde{\alpha }-\alpha \Vert \leq d_{\alpha }}\, \left|f(\tilde{\alpha })-l(\tilde{\alpha })\right|\, <\, \varepsilon \, ,\]
holds, where $\, \inf \, $ is taken over all linear functions $l$.
Then\[
f\left(\alpha +M(\alpha )\right)-f(\alpha )<-5\varepsilon \qquad \Longrightarrow \qquad \sE \left(\, f(\xi _{n+1})-f(\alpha )\, |\, \xi _{n}=\alpha \right)<-\varepsilon \, .\]

\end{lemma}

\end{document}